\documentclass[11pt,reqno,oneside]{amsart}

\newtheorem{defin}{Definition}[section]
\newtheorem{theorem}[defin]{Theorem} 
\newtheorem{lemma}[defin]{Lemma} 

\def\C{{\mathbb  C}}
\def\N{{\mathbb  N}}
\def\R{{\mathbb  R}}
\def\B{{\mathbb  B}}

\def\cJ{{\mathcal J}}
\def\cC{{\mathcal C}}

\def\d{\delta}
\def\e{\varepsilon}
\def\z{\zeta}

\def\l{\lambda}

\def\bar{\overline}        

\def\disc{\triangle}             
\def\cdisc{\bar\triangle}        
\def\bdisc{b\disc}               
\def\cn{{\mathbb C}^n}

\def\ctwo{{\mathbb C}^2}

\def\dim{{\rm dim}\,}    
\def\st{such that}

\def\psh{plurisubharmonic}               

\def\holo{holomorphic}
\def\nbd{neighborhood}
\def\dist{{\rm dist}}

\def\ss{\subset\!\subset}   

\begin{document}
\bibliographystyle{abbrv}
\title{\bf On proper discs in complex manifolds}         

\author{Barbara Drinovec Drnov\v sek}

\address{Institute of Mathematics, Physics and Mechanics, University of Ljubljana,
Jadranska 19, SI-1000 Ljubljana, Slovenia}
\email{Barbara.Drinovec@fmf.uni-lj.si}
\date{\today}

\begin{abstract}
Let $X$ be a complex manifold of dimension at least $2$ which has an 
exhaustion function whose Levi form has at each point at least $2$ 
strictly positive eigenvalues.
We construct proper holomorphic discs in $X$
through any given point and in any given direction.
\end{abstract}

\maketitle

\section{Introduction and the results}

Denote by $\disc$ the open unit disc in $\C$. 
Our main result is the following

\begin{theorem}
Let $X$ be a complex manifold of dimension $n\ge 2$ which has an 
exhaustion function whose Levi form has at each point at least $2$ 
positive eigenvalues.
Given $p\in X$ and a vector $v$ tangent to $X$ at $p$,
there is a proper holomorphic map
$f\colon\disc\to X$ such that $f(0)=p$ and $f'(0)=\l v$ for some $\l>0$.
\label{izrek0}
\end{theorem}

Note that in the case $n=2$ the manifold $X$ having the property in the theorem is Stein.
It is known that if $X$ is a Stein manifold  of dimension at least $2$,
then for each point $p$ in  $X$ and for each vector $v$ tangent to $X$ at $p$,
there is a proper holomorphic map $f\colon\triangle\to X$ such that $f(0)=p$ 
and $f'(0)=\l v$ for some $\l>0$ \cite{Glo,FG1}.
Therefore, our result is new in the case $n>2$.

In the theory of $q$-convex manifolds the manifolds with the above property are 
called $(n-1)$-complete manifolds. For general theory of $q$-convex manifolds 
we refer to \cite{Gra,HL}. The survey \cite{Col} gives a list of examples 
(with references) and some open questions concerning $q$-convexity.

The conclusion of the theorem is not valid in general for complex manifolds 
of dimension at least $2$ which have an 
exhaustion function whose Levi form has at each point at least $1$ 
positive eigenvalue.
Indeed, in \cite{FG1} the authors constructed for every $n\ge 2$ 
a smoothly bounded domain $\Omega\ss\cn$ and a point $p\in \Omega$ \st\ 
there is no proper holomorphic map $f\colon\disc\to \Omega$ with $p\in f(\disc)$. 
Since every noncompact complex manifold of dimension $n\ge 2$  has an 
exhaustion function whose Levi form has at each point at least $1$ 
positive eigenvalue \cite{GW}
we get the desired example for which the conclusion of Theorem \ref{izrek0} fails.

A. Dor \cite{Dor} proved that there exists a bounded domain $\Omega $ in 
$\C^n$ $(n \ge 2)$ such that there is no proper holomorphic 
mapping from the unit disc to $\Omega $. 

We prove the following result.

\begin{theorem}
Let $X$ be a complex manifold of dimension $n\ge 2$ and
assume that $\rho\colon X\to\R$ is an exhaustion function
whose Levi form has at least $2$ 
positive eigenvalues at each point of the set $\{\rho>M\}$ for some $M\in \R$.
Let $d$ be a complete metric on $X$ which induces the manifold topology.
Given $\e>0$, $0<r<1$, and a continuous map $f\colon\cdisc\to X$ such that 
$\rho(f(\z))>M$ $(\z\in\bdisc)$
there is a proper holomorphic map
$g\colon\disc\to X$ such that
\item[(i)] $d(g(\z),f(\z))<\e$ for $|\z|<r$,
\item[(ii)] $g(0)=f(0)$,
\item[(iii)] $g'(0)=\l f'(0)$ for some $\l>0$.
\label{izrek1}
\end{theorem}

Since there are small holomorphic discs through any given point in any given direction
on each complex manifold, Theorem \ref{izrek1} easily implies Theorem \ref{izrek0}.

If $\dim X=2$ then the function $\rho$ is
strictly \psh\ in the set $\{\rho >M\}$ and in this case theorem was proved in
\cite{Glo},\cite[Theorem 1.1]{FG2}. So we only need to treat the case $\dim X\ge 3$.

In the proof of Theorem \ref{izrek1} we shall push the boundary of a given analytic disc 
outside a given sublevel set of $\rho$. 
Since our manifold does not necessarily lie in Euclidean space, we are not
able to do this by adding a suitable polynomial map as it was done in \cite{Glo}. 
Instead we use convex bumps. At the first step 
we push the boundary outside the given sublevel set union one bump. 
At the next step the boundary lies outside the sublevel set union two bumps.
These bumps are constructed in such a way that they fill the space 
between two level sets of $\rho$.
In a finite number of steps the boundary of the disc lies outside 
the bigger sublevel set.

In section \ref{sec2} we prove that we can push
the boundary of a given holomorphic disc along a continuous family
of small holomorphic discs attached to the boundary;
these small discs are not constant only in a fixed coordinate \nbd . 
To do this we first solve approximately a Riemann-Hilbert boundary value problem 
\cite[Lemma 5.1]{Glo}, \cite{FG1}  to get a holomorphic
map from the part of the open unit disc which was initially mapped into the fixed
coordinate \nbd\ and then we obtain the new holomorphic disc as a 
solution of a nonlinear Cousin problem (due to
J.-P. Rosay \cite{Ros}).

In section \ref{sec3} we construct convex bumps which provide
the continuous family of holomorphic discs, and for this family we use the result 
from section \ref{sec2}. We prove Theorem \ref{izrek1} in section \ref{sec4}.

\section{The main lemma}
\label{sec2}

We will need the following elementary lemma
\begin{lemma}
Let $M$ be a metric space with distance function $d$.
Let $U_1,U_2$ be open sets in $\C$ such that $\cdisc\subset U_1\cup U_2$
and assume that $f_k\colon \bar{U_k\cap\disc}\to M$ $(k=1,2)$ are continuous maps.
Further, assume that for some $\e>0$ there is a continuous map $g_0\colon\disc\to M$
such that
$$d(g_0(\z),f_k(\z))<\e \ (\z\in U_k\cap\disc, k=1,2).$$
Then there is $R$, $0<R<1$, so close to $1$ that
$$d(g_0(R\z),f_k(\z))<5\e\ (\z\in U_k\cap\cdisc, k=1,2).$$
\label{element}
\end{lemma}

\begin{proof}
It is easy to see that $(U_1\setminus U_2)\cap \disc\ss U_1$ and $(U_2\setminus U_1)\cap \disc\ss U_2$. Therefore
there is $R$, $0<R<1$, so close to $1$ that 
\begin{equation}
R((U_1\setminus U_2)\cap \disc)\subset U_1\text{ and } 
R((U_2\setminus U_1)\cap \disc)\subset U_2.
\label{elem1}
\end{equation}
Since $f_k$ is uniformly continuous on $\bar{U_k\cap\disc}$ $(k=1,2)$ by increasing $R$, $R<1$, we obtain
\begin{equation}
\text{ if }\z,\eta\in\bar{U_k\cap \disc},\ |\z-\eta|\le 1-R\text{ then } d(f_k(\z),f_k(\eta))<\e
\ (k=1,2).
\label{elem2}
\end{equation}
For $\z\in(U_1\setminus U_2)\cap\disc$ using (\ref{elem1}) and (\ref{elem2})
we get
$$d(g_0(R\z),f_1(\z))\le d(g_0(R\z),f_1(R\z))+d(f_1(R\z),f_1(\z))<2\e,$$
and similarly for $\z\in(U_2\setminus U_1)\cap\disc$.
Take $\z\in U_1\cap U_2\cap\disc$. 
If $R\z\in U_1$ then as above we get $d(g_0(R\z),f_1(\z))<2\e$.
If $R\z\notin U_1$ then $R\z\in U_2$,
and there is $t$, $R<t<1$, such that $t\z\in U_1\cap U_2$, and we have
$ d(g_0(R\z),f_1(\z))\le d(g_0(R\z),f_2(R\z))+d(f_2(R\z),f_2(t\z))+
 d(f_2(t\z),g_0(t\z))+d(g_0(t\z),f_1(t\z))+d(f_1(t\z),f_1(\z))<5\e.$
Similarly $d(g_0(R\z),f_2(\z))<5\e$.
This completes the proof.
\end{proof}

The following lemma, which holds for any complex manifold of dimension at least two,
is the main tool in the inductive construction of a proper holomorphic disc.
The proof depends on the solution of a nonlinear Cousin problem due to
J.-P. Rosay \cite{Ros}.

\begin{lemma}
Let $X$ be a complex manifold of dimension $n\ge 2$ endowed with a Riemannian metric, 
which induces the distance function $d$ on $X$.
Assume that $\Omega\ss X$ is an open coordinate \nbd\ and $\Omega_0\ss \Omega$. 
Let $f$ be a holomorphic map from a \nbd\ of $\cdisc$ to $X$.
Let $K$ be a compact subset of $X$ and $r$, $0<r<1$, \st\ $K\cap \Omega=\emptyset$ and
$K\cap f(\cdisc\setminus r\disc)=\emptyset$.
Assume that $H\colon\bdisc\times\cdisc\to X$ is a continuous map with the following 
properties:
\item[(i)] for each $\z\in\bdisc$ the map $\eta\mapsto H(\z,\eta)$ is holomorphic on $\disc$,
\item[(ii)] $H(\zeta,0)=f(\zeta)$ $(\zeta\in\bdisc)$, 
\item[(iii)] if for some $\z\in\bdisc$ we have $f(\z)\notin \Omega_0$ then 
$H(\z,\eta)=f(\zeta)$ $(\eta\in\cdisc)$,
\item[(iv)] if for some $\z\in\bdisc$ we have $f(\z)\in \Omega_0$ then 
$H(\z,\eta)\in \Omega$ $(\eta\in\cdisc)$.

Given $\e>0$, there is a 
continuous map $g\colon\cdisc\to X$, \holo\ on $\disc$, \st\ 
\item[(i')] $d(g(\z),H(\z,\bdisc))<\e$ $(\z\in\bdisc)$,
\item[(ii')] $d(g(\z),f(\z))<\e$ $(\z\in r\cdisc)$,
\item[(iii')] $g(\cdisc\setminus r\disc)\cap K=\emptyset$,
\item[(iv')] $g(0)=f(0)$,
\item[(v')] $g'(0)=Rf'(0)$ for some $R$, $r<R<1$.
 \label{lift}
\end{lemma}

\begin{proof}
By assumption there is $\rho>1$ such that
$f$ is \holo\ on $\rho\disc$. Choose $\rho'$, $1<\rho'<\rho$. Define the map 
$\tilde f\colon\rho\disc\to X\times\C$ by $\tilde f(\z)=(f(\z),\z)$. 
The map $\tilde f$ is a holomorphic embedding (not proper), so there is an open
\nbd\ $\tilde \Omega_1$ of $\tilde f(\rho'\disc)$ in $X\times \C$
and a bi\holo\ map $\tilde\Phi_1$ from $\tilde \Omega_1$
onto a bounded open subset of $\C^{n+1}$ (see \cite{Roy}, \cite{LS}, \cite[Lemma 1.1]{Ro2}).
With no loss of generality one can in addition assume that
the derivative  of $\tilde\Phi_1$ at $0$ is the identity map.
Choose a compact \nbd\ $\tilde K_1$ of $\tilde f(\cdisc)$ in $\tilde \Omega_1$.

Choose a biholomorphic map $\Phi_2$ from $\Omega$ to an open subset of $\cn$ and
an open set $\Omega_2$ such that $\Omega_0\ss \Omega_2\ss \Omega$, and
$H(f^{-1}(\bar \Omega_0),\cdisc)\subset \Omega_2$. 
Let $\tilde \Omega_2=\Omega\times 3\disc$ and $\tilde K_2=\bar\Omega_2\times 2\cdisc$.
Let $\tilde \Phi_2(z,\z)=(\Phi_2(z),\z)$ for $(z,\z)\in \tilde \Omega_2$ and note that
$\tilde \Phi_2$ maps $\tilde \Omega_2$ biholomorphically into $\C^{n+1}$.

By decreasing $\e>0$ if necessary we may assume that 
\begin{eqnarray}
d(K,f(\cdisc\setminus r\disc))>\e \text{ and }
d(K,\Omega_2)>\e .
\label{1b}
\end{eqnarray}

Denote by $\B$ the open unit ball in $\cn$
and by $\Pi\colon X\times \C\to X$ the canonical projection to the first factor.
There is $\alpha>0$ so small that
 \begin{eqnarray}
 \Phi_2(H(f^{-1}(\bar \Omega_0),\cdisc))+\alpha\B\subset \Phi_2(\Omega_2),\label{eta1}\\
  \left.\begin{array}{c}
 \text{ if } z\in \tilde\Phi_1(\tilde K_1),\ z'\in\C^{n+1},\ 
 |z-z'|<\alpha, \\
 \text{ then } d(\Pi(\tilde\Phi_1^{-1}(z)),\Pi(\tilde\Phi_1^{-1}(z')))<\tfrac{\e}{6},
 \end{array}\right\} \label{eta}\\
 \left.\begin{array}{c}
\text{ if } z\in \tilde\Phi_2(\tilde K_2),\ z'\in\C^{n+1},\ 
 |z-z'|<\alpha,\\ \text{ then }
 d(\Pi(\tilde\Phi_2^{-1}(z)),\Pi(\tilde\Phi_2^{-1}(z')))<\tfrac{\e}{6}.
 \end{array}\right\}
 \label{eta222}
 \end{eqnarray}

By slightly enlarging $\Omega_0$ we may assume that
either $f(\bdisc)\subset \bar \Omega_0$ or
the set $f^{-1}(\bar \Omega_0)\cap\bdisc$ is at most finite union of disjoint closed arcs.
Consider first the second case: 
denote these arcs by $\{I_j\}_{j\in \cJ}$ where $\cJ$ is finite and
where $I_j$ are pairwise disjoint. 
For each $j\in\cJ$ one can find a smooth simple closed curve 
$\Gamma_j\subset \cdisc\setminus r\cdisc$ such that $\Gamma_j\cap\bdisc$ is a \nbd\ of $I_j$
in $\bdisc$ and $\Gamma_j$ are pairwise disjoint. 
Each $\Gamma_j$ bounds a domain $D_j\subset\disc\setminus r\cdisc$, which is
conformally equivalent to the unit disc. 
Since $f(I_j)\subset \Omega_0$ one can choose $\Gamma_j$ in such a way that,
in addition to the above, we have
$f(\bar D_j)\subset \Omega_2$ $(j\in \cJ)$.
Choose  a homeomorphic map $h_j$
from $\cdisc$ to $\bar D_j$, which is holomorphic on $\disc$, 
and let $V_j=h_j(\{th_j^{-1}(\z);\z\in I_j, t\in[0,1]\})$.
One can choose an  open neighborhood $W_j$ of $V_j$ in $\C$
such that $W_j\cap \disc\subset D_j$.
Denote by $U_1$ the set $\rho'\disc\setminus \bar{\cup_j V_j}$ and 
by $U_2$ the union $\cup_j W_j$.

If $f(\bdisc)\subset \bar \Omega_0$ then let $I_1=\bdisc$ and choose $r_0$, $r<r_0<1$,
such that $f(\cdisc\setminus r_0\disc)\subset \Omega_2$. Let $U_1=\frac{1+r_0}{2}\disc$
and let $U_2=W_1=\rho'\disc\setminus \bar{r_0\disc}$.

Now we are in the situation of Section 5 in \cite{Ros}: 
The sets $U_1$ and $U_2$ are open in $\C$ and satisfy $\cdisc\subset U_1\cup U_2$.
Let $\omega_{12}'=\tilde\Phi_1(\tilde\Omega_1\cap\tilde\Omega_2)$ and let
$\omega_{12}$ be the image under $\tilde\Phi_1$ of a \nbd\ of 
$\tilde K_1\cap\tilde K_2$ such that $\omega_{12}\ss\omega_{12}'$. 
By \cite[Proposition 1']{Ros}  there is a $\delta>0$
with the following property: if $u_k$ is a holomorphic 
map from $U_k\cap \disc$ into $\C^{n+1}$ $(k\in\{1,2\})$
such that $u_1(U_1\cap U_2\cap\disc)\subset \omega_{12}$, and
$|u_2(\z)-(\tilde\Phi_2\circ\tilde\Phi_1^{-1}\circ u_1)(\z)|\le \d$  $(\z\in U_1\cap U_2\cap \disc)$,
then there are holomorphic maps
$v_k$ from $U_k\cap\disc$ into $\C^{n+1}$ such that
$|v_k(\z)|\le \alpha$ $(\z\in U_k\cap\disc)$ for $k=1,2$ and 
$u_2+v_2=\tilde\Phi_2\circ\tilde\Phi_1^{-1}\circ (u_1+v_1)$
on $U_1\cap U_2\cap\disc$.
Moreover, one can impose $v_1(0)=0$.
By a minor change in the proof of \cite[Proposition 1']{Ros} one can further impose 
$v_1'(0)=0$. For the sake of completeness we provide the details.
We will adapt the same notations as in the proof of \cite[Proposition 1']{Ros}.
We only need to change the solution operator $T$ which solves the standard 
additive Cousin problem.
This continuous linear operator associates to a bounded holomorphic map 
$\alpha_{12}\in (H^{\infty}(U_1\cap U_2\cap\disc))^{n+1}$ holomorphic maps
$T_j(\alpha_{12})\in (H^{\infty}(U_j\cap\disc))^{n+1}$ $(j=1,2)$ such that
$\alpha_{12}=T_1(\alpha_{12})-T_2(\alpha_{12})$.
We need that in addition to the above, it satisfies
$T_1(\alpha_{12})(0)=0$ and $T_1(\alpha_{12})'(0)=0$.
This property of operator $T$ implies that $v_1'(0)=0$.
Following the proof of \cite[Theorem 1.4.5]{Hor} the solution of the additive Cousin 
problem is reduced to solving $\bar\partial$-equation as follows.
We can choose a cut off function $\varphi$ such that ${\rm supp}\, \varphi\ss U_1$ and
$\varphi$ equals $1$ on $(U_1\setminus U_2)\cap\cdisc$. The solution of the additive Cousin problem
is of the form
\begin{equation}
T_1(\alpha_{12})=(1-\varphi)\alpha_{12}+u \ ,
\ T_2(\alpha_{12})=-\varphi\alpha_{12}+u \ , 
\label{T1T2}
\end{equation}
where $u$ is a solution of $\bar\partial$-equation to assure that the maps $T_j(\alpha_{12})$
 are holomorphic on $U_j\cap \disc$ $(j=1,2)$. 

Instead of solving (\ref{T1T2}) we will solve the following
\begin{eqnarray*}
T_1(\alpha_{12})(\z)=&(1-\varphi(\z))\alpha_{12}(\z)+\z^2u(\z)\ &(\z\in U_1\cap\disc),\\
T_2(\alpha_{12})(\z)=&-\varphi(\z)\alpha_{12}(\z)+\z^2u(\z) \ &(\z\in U_2\cap\disc).
\end{eqnarray*}
Therefore, if the map $u$ satisfies
$$\frac{\partial u}{\partial \bar\z}(\z)=
\dfrac{\alpha_{12}(\z)}{\z^2}\frac{\partial \varphi}{\partial \bar\z}(\z)
\ \ (\z\in\disc),$$
then the maps $T_j(\alpha_{12})$ $(j=1,2)$ are holomorphic.
By the properties of $\varphi$ one can solve the above
$\bar\partial$-equation with estimates to prove that $T$ is continuous and linear.

Now we proceed with the proof of Lemma \ref{lift}.
Let $u_1(\z)=\tilde\Phi_1(f(\z),\z)$ $(\z\in U_1\cap \cdisc)$.
In the next paragraph we shall define map $u_2$.

First consider the case that 
$f^{-1}(\bar \Omega_0)\cap\bdisc$ consists of finitely many closed arcs denoted
by $I_j$.
For each $j\in \cJ$ define the map $H_j\colon\Gamma_j\times\cdisc\to X$ by
$H_j(\z,\eta)=H(\z,\eta)$ $(\z\in I_j,\, \eta\in\cdisc)$ and $H_j(\z,\eta)=f(\z)$ 
$(\z\in\Gamma_j\setminus I_j,\, \eta\in\cdisc)$. 
Property (iii) implies that $H_j$ is continuous, and by (i),
$H_j$ is holomorphic in the second variable.
By (\ref{eta1}), and as $f(\bar D_j)\subset \Omega_2$, there is an $\alpha_j$, $0<\alpha_j<\alpha$,
so small that
\begin{equation}
\Phi_2(H_j(\Gamma_j,\cdisc))+\alpha_j\B\subset \Phi_2(\Omega_2).
\label{etaj}
\end{equation}
Let $\gamma>0$.
Using
\cite[Lemma 5.1]{Glo} for the map $G_j(\z,\eta)=\Phi_2(H_j(h_j(\z),\eta))-\Phi_2(f(h_j(\z)))$ 
$((\z,\eta)\in\bdisc\times\cdisc)$
we get the polynomial map $P_j$ from $\C$ into $\cn$ with the following properties
\begin{eqnarray*}
 P_j(\z)&\in& G_j(\z,\bdisc)+\gamma\B,\ (\z\in \bdisc),\\
 P_j(t\z)&\in& G_j(\z,\cdisc)+\gamma\B,\ (\z\in \bdisc, 0\le t\le 1).
\end{eqnarray*}
If $\gamma>0$ is small enough then the map $p_j\colon \bar D_j\to \cn$,
defined by $p_j(\z)=P_j(h_j^{-1}(\z))+\Phi_2(f(\z))$ $(\z\in\bar D_j)$,
is continuous, holomorphic on $D_j$, and satisfies the following
\begin{eqnarray}
p_j(\z)\in \Phi_2(H_j(\z,\bdisc))+{\alpha}\B \ (\z\in \Gamma_j),\label{en00}\\ 
p_j(\z)\in \Phi_2(H_j(\Gamma_j,\cdisc))+\alpha_j\B\  (\z\in \bar D_j), \label{en0}\\ 
|p_j(\z)-\Phi_2(f(\z))|<{\d} \ (\z\in D_j\setminus V_j).\label{en1}
\end{eqnarray}
Note that by (\ref{etaj}) and (\ref{en0}), we get $p_j(\bar D_j)\subset \Phi_2(\Omega_2)$.

Define the map $u_2\colon U_2\cap \cdisc\to \C^{n+1}$ by 
$u_2(\z)=(p_j(\z),\z)$ $(\z\in \bar{W_j\cap\disc})$.

In the simpler case when $f(\bdisc)\subset \bar \Omega_0$ we solve
approximately the Riemann-Hilbert problem for the map 
$G_1(\zeta,\eta)=\Phi_2(H(\zeta,\eta))-\Phi_2(f(\z))$ $((\zeta,\eta)\in\bdisc\times\cdisc)$
such that the solution polynomial map $P_1$ as above satisfies
\begin{equation}
|P_1(\z)|<\min\{\delta,\alpha\} \text{ for } \ |\z|<\frac{1+r_0}{2}.\label{primer}
\end{equation}
Similarly as above we define
$p_1(\z)=P_1(\z)+\Phi_2(f(\z))$ $(\z\in U_2\cap \cdisc)$ and
$u_2(\z)=(p_1(\z),\z)$ $(\z\in U_2\cap \cdisc)$.

Property (\ref{en1}) (or (\ref{primer})) implies
that $|u_2(\z)-(\tilde\Phi_2\circ\tilde\Phi_1^{-1}\circ u_1)(\z)|\le \d$  
$(\z\in U_1\cap U_2\cap \disc)$.
Therefore by the above there exist holomorphic maps $v_k\colon U_k\cap\disc\to \C^{n+1}$
$(k=1,2)$ such that 
\begin{eqnarray*}
|v_k(\z)|<\alpha \ (\z\in U_k\cap\disc,\ k=1,2),\\
u_2+v_2=\tilde\Phi_2\circ\tilde\Phi_1^{-1}\circ (u_1+v_1) \text{ on } U_1\cap U_2\cap\disc,\\
v_1(0)=0,\ v_1'(0)=0.
\end{eqnarray*}
This implies that one can define 
$$g_0(\z)=\Pi(\tilde\Phi_k^{-1}((u_k+v_k)(\z)))\ (\z\in U_k\cap\disc,k=1,2).$$
The map $g_0$ 
is holomorphic on $\disc$, $g_0(0)=f(0)$, $g_0'(0)=f'(0)$, 
and by (\ref{eta}) and (\ref{eta222}), we have
\begin{eqnarray}
d(g_0(\z),f(\z))<\tfrac{\e}{6}\ (\z\in U_1\cap \disc),\label{raz1}\\
d(g_0(\z),\Phi_2^{-1}(p_j(\z)))<\tfrac{\e}{6}\ (\z\in W_j\cap\disc,\ j\in\cJ).\label{raz2}
\end{eqnarray}

By Lemma \ref{element}, there is $R$, $r<R<1$, so close to $1$ that the map
$g(\z)=g_0(R\z)$ $(\z\in\cdisc)$ satisfies the following
\begin{eqnarray}
d(g(\z),f(\z))<\tfrac{5\e}{6} \ (\z\in U_1\cap \cdisc),\label{raz3}\\
d(g(\z),\Phi_2^{-1}(p_j(\z)))<\tfrac{5\e}{6} \ (\z\in W_j\cap\cdisc,\ j\in\cJ).\label{raz4}
\end{eqnarray}
The map $g$ is continuous on $\cdisc$ and 
holomorphic on $\disc$. 
To prove (i') choose $\z\in \bdisc\cap U_1$.
Since $f(\bdisc \cap U_1)\cap\Omega_0=\emptyset$, properties
(iii) and (\ref{raz3}) imply that
$d(g(\z),H(\z,\bdisc))=d(g(\z),f(\z))<\e.$
For  $\z\in U_2\cap \bdisc$ it holds that $\z\in W_j$ for some $j\in\cJ$ and by 
(\ref{raz4}), (\ref{en00}) and (\ref{eta222}) we get
$$d(g(\z),H(\z,\bdisc))\le d(g(\z),\Phi_2^{-1}(p_j(\z)))+
d(\Phi_2^{-1}(p_j(\z)),H(\z,\bdisc))<\e.$$
So (i') holds.
Since $r\cdisc\subset U_1$  property (\ref{raz3}) implies (ii').
Take $\z\in\cdisc\setminus r\disc$. If $\z\in U_1$, then (\ref{raz3}) and (\ref{1b})
imply that $g(\z)\notin K$. If $\z\in U_2$, then  $\z\in W_j$ for some $j\in\cJ$ and
since $p_j(\z)\in \Phi_2(\Omega_2)$ properties (\ref{raz4}) and (\ref{1b}) imply that 
$g(\z)\notin K$. Thus we proved (iii').
It is easy to see that $g(0)=f(0)$ and $g'(0)=Rf'(0)$, which proves (iv') and (v').
This proves the lemma.
\end{proof}

\section{Convex bumps}
\label{sec3}

We will construct convex bumps introduced by Grauert (see \cite{HL}) in order 
to provide a continuous family of analytic discs,
which is needed to use Lemma \ref{lift} in the inductive proof of Theorem \ref{izrek1}.

We denote by $d_{1,2}$ the partial differential with respect to
the first two complex coordinates on $\cn$.

Let $X$ be a $n$-dimensional complex manifold and 
let $A,B\subset X$ be relatively compact open
sets in $X$. We say that {\it $B$ is a $2$-bump on $A$} if there exist an open 
set $\Omega\subset X$ containing $\bar B$, a biholomorphic map $\Phi$ from $\Omega$ onto a 
convex subset $\omega$ in $\cn$, and smooth real functions $\rho_B\le \rho_A$ on $\omega$  such that
$$\Phi(A\cap \Omega)=\{z\in \omega;\, \rho_A(z)<0\},\ 
\Phi((A\cup B)\cap \Omega)=\{z\in \omega;\, \rho_B(z)<0\},$$
and the functions $\rho_A$ and $\rho_B$ are strictly plurisubharmonic with respect to the first two coordinates.

We say that {\it $B$ is a convex $2$-bump on $A$} if, in addition to the above, 
$\rho_A$ and $\rho_B$ are strictly convex in $z_1,z_2$ (with respect to the underlying real
coordinates), and $d_{1,2}(t\rho_A+(1-t)\rho_B)$ 
is non degenerate on $\omega$ for each $t\in[0,1]$.

\begin{lemma}
Let $X$ be a complex manifold of dimension $n\ge 2$, equipped with some metric $d$.
Let $A,B\subset X$ such that $B$ is a convex $2$-bump on $A$ and let $K$ be a compact 
subset of $A$.
Assume that $f\colon\cdisc\to X$ is a continuous map, \holo\ on $\disc$, such that
$f(\bdisc)\cap \bar A=\emptyset$, and $r$, $0<r<1$, such that 
$f(\cdisc\setminus r\disc) \cap K=\emptyset$.
Given $\e>0$ there is a continuous map $g\colon\cdisc\to X$, \holo\ on $\disc$,
with the following properties
\item[(i)] $g(\bdisc)\cap (\bar{A\cup B})=\emptyset$,
\item[(ii)] $d(g(\z),f(\z))<\e$ $(|\z|\le r)$,
\item[(iii)] $g(\cdisc\setminus r\disc)\cap K=\emptyset$,
\item[(iv)] $g(0)=f(0)$, 
\item[(v)] $g'(0)=Rf'(0)$ for some $R$, $r<R<1$.
\label{K1}
\end{lemma}

\begin{proof} We may assume that $f$ is holomorphic in a \nbd\ of $\cdisc$.
Since $B$ is a convex $2$-bump on $A$ there are
a biholomorphic map $\Phi\colon \Omega\to \cn$ onto a convex subset $\omega$ of $\cn$,
and smooth functions $\rho_A,\rho_B\colon \omega\to \R$ such that
$$\Phi(A\cap \Omega)=\{z\in \omega;\, \rho_A(z)<0\},\ 
\Phi((A\cup B)\cap \Omega)=\{z\in \omega;\, \rho_B(z)<0\},$$
and the functions  $\rho_A$ and $\rho_B$ are strictly convex with respect to the first 
two coordinates, and $d_{1,2}(t\rho_A+(1-t)\rho_B)$ 
is non degenerate on $\omega$ for each $t\in[0,1]$.  

We can choose $\l>0$ so small that
$\Phi(f(\bdisc)\cap \Omega)\cap \{z\in \omega;\rho_A(z)\le\l \}=\emptyset$, the set
$\omega_0=\{z\in \omega;\rho_A(z)>\l,\rho_B(z)\le\l\}$ is relatively compact in $\omega$, 
and $\bdisc\cap f^{-1}(\Phi^{-1}(\{z\in \omega;\, \rho_B(z)\le \l\}))$ is either $\bdisc=:I_1$ or a union of finitely many closed arcs which we denote by $I_j$.

Choose a point $q\in \omega_0$ and
write $q=(q_1,q_2,q'')$.
There is exactly one $\mu$, $0\le \mu\le 1$, such that 
$\mu(\rho_A(q)-\l)+(1-\mu)(\rho_B(q)-\l)=0$. 
The function $\mu(\rho_A-\l)+(1-\mu)(\rho_B-\l)$ is defined on $\omega$ and 
it is strictly convex in the first two coordinates.
Denote by $M_{\mu,q''}$ the set $\{(z_1,z_2,q'')\in \omega;
\mu(\rho_A(z_1,z_2,q'')-\l)+(1-\mu)(\rho_B(z_1,z_2,q'')-\l)=0\}$.
Note that $M_{\mu,q''}$ is a real submanifold of dimension $3$ in $\ctwo\times\{q''\}$.
Denote by $T_qM_{\mu,q''}$ its real tangent space at $q$.
The intersection $E_q=T_qM_{\mu,q''}\cap\imath T_qM_{\mu,q''}$ is a complex line.
By strict convexity the intersection of $\{q\}+E_q$ with
$\{z\in \Omega;\rho_B(z)\le \l\}$ is a bounded connected convex subset of $\{q\}+E_q$
therefore it is conformally equivalent to the unit disc. 
If we vary $q$ smoothly these convex sets
vary smoothly.
The set $\Phi(f( I_j))$ is contained in $\omega_0$ for each $j$.
Therefore with a proof similar to the proof of \cite[Lemma 4.1]{Glo} we obtain 
a continuous map $H_j\colon  I_j\times\cdisc\to \cn$ such that
\item[(a)] for each $\z\in I_j$ the map $H_j(\z,\eta)$ is holomorphic in $\eta$,
\item[(b)] $H_j(\z,0)=\Phi(f(\z))$ $(\z\in  I_j)$,
\item[(c)] $H_j( I_j,\cdisc)\subset \omega$,
\item[(d)] $\rho_B(H_j(\z,\eta))=\l$ $(\z\in  I_j,\eta\in b\disc)$,
\item[(e)] if $\rho_B(\Phi(f(\z)))=\l$ then $H_j(\z,\eta)=\Phi(f(\z))$ 
$(\z\in I_j,\eta\in\cdisc)$.

We define a map $H\colon\bdisc\times\disc\to X$ by
\begin{equation*}
H(\z,\eta)=\left\{\begin{array}{ll}\Phi^{-1}(H_j(\z,\eta)),&\z\in I_j,\\
f(\z),& \z\in\bdisc\setminus \cup_j  I_j.\end{array}\right. 
\end{equation*}
The map $H$ is continuous by construction and it satisfies the following
\item[(a')] for each $\z\in\bdisc$ the map $H(\z,\eta)$ is holomorphic in $\eta$,
\item[(b')] $H(\zeta,0)=f(\zeta)$ $(\zeta\in\bdisc)$, 
\item[(c')] if for some $\z\in\bdisc$ we have $f(\z)\notin \Phi^{-1}(\omega_0)$ then 
$H(\z,\eta)=f(\z)$ $(\eta\in\cdisc)$,
\item[(d')] if for some $\z\in\bdisc$ we have $f(\z)\in \Phi^{-1}(\omega_0)$ then 
$H(\z,\cdisc)\subset \Phi^{-1}(\omega)$ and $H(\z,\eta)\in 
\Phi^{-1}(\{z\in \omega;\rho_B(z)=\l,\rho_A(z)\ge\l\})$ $(\eta\in\bdisc)$.

Take $\e_0$, $0<\e_0<\e$, so small that
\begin{eqnarray}
d(\Phi^{-1}(\{z\in \omega;\rho_B(z)=\l,\rho_A(z)\ge\l\}),\Phi^{-1}(\{z\in \omega;\rho_B(z)\le 0\}))>\e_0 .
\label{epsilon0}
\end{eqnarray}

Now we use Lemma \ref{lift} to get the map $g$ such that
\item[(i')] $d(g(\z),H(\z,\bdisc))<\e_0$ $(\z\in\bdisc)$,
\item[(ii')] $d(g(\z),f(\z))<\e_0$ $(\z\in r\cdisc)$,
\item[(iii')] $g(\cdisc\setminus r\disc)\cap K=\emptyset$,
\item[(iv')] $g(0)=f(0)$,
\item[(v')] $g'(0)=Rf'(0)$ for some $R$, $r<R<1$.

Properties (ii)-(v) follow from (ii')-(v').
By (\ref{epsilon0}), (d') and (i') we get (i).
\end{proof}

\section{Proof of Theorem \ref{izrek1}}
\label{sec4}

As we have already explained in the introduction
we only need to treat the case $\dim X\ge 3$.

\begin{lemma}
Let $X$ be a complex manifold of dimension $n\ge 3$.
Let $\Omega\ss X$ and let $\rho\colon \Omega\to \R$ be a smooth function
such that $\{z\in \Omega;\, a\le \rho(z)\le b\}\ss \Omega$ and
such that  the Levi form of $\rho$ has at each point at least $2$ positive eigenvalues.
Assume that $\rho$ has at most one critical point  in $\{z\in \Omega;\, a\le \rho(z)\le b\}$
and, if $q$ is a critical point of $\rho$, then further assume that $a< \rho(q)< b$, and
that $q$ is a non-degenerate critical point.
Let $K$ be a compact subset of $X$ \st\ $K\cap \Omega=\emptyset$. 
Assume that $f\colon\cdisc \to X$ is a 
continuous map, \holo\ on $\disc$, \st\ $f(\bdisc)\subset \{z\in \Omega;\rho(z)>a\}$
and choose $r$, $0<r<1$, such that $f(\cdisc\setminus r\disc)\cap K=\emptyset$.
Given $\e>0$ there exists a continuous map $g\colon\cdisc\to X$, \holo\ on $\disc$,
with the following properties 
\item[(i)] $g(\bdisc)\subset \{z\in \Omega;\rho(z)>b\}$,
\item[(ii)] $d(g(\z),f(\z))<\e$ $(|\z|\le r)$,
\item[(iii)] $g(\cdisc\setminus r\disc)\cap K=\emptyset$,
\item[(iv)] $g(0)=f(0)$,
\item[(v)] $g'(0)=\l f'(0)$ for some $\l$, $r<\l<1$.
\label{kritdvig}
\end{lemma}

\begin{proof}
Note that, if $\rho$ is a smooth function defined on a complex manifold $X$,
whose Levi form has at least $2$ positive eigenvalues 
at some point $w\in X$, then there are holomorphic coordinates near $w$ such that $\rho$ is
strictly plurisubharmonic in $z_1$,$z_2$. Moreover, if $w$ is a regular point of $\rho$, then
Narasimhan's lemma on local convexification implies that, in local holomorphic coordinates,
$\rho$ can be made strictly convex in $z_1$,$z_2$.
Both conditions are stable under small perturbations. 
If $\rho$ does not have local minima in $\Omega$, then by
\cite[Lemma 12.3]{HL} we get finitely many domains 
$\{z\in \Omega; \rho(z)<a\}=A_0\subset A_1\subset\cdots\subset A_m=\{z\in \Omega; \rho(z)<b\}$
in $X$ such that for every $k=0,1,\ldots,m-1$ we have 
$A_{k+1}=A_k\cup B_k$, where $B_k$ is a $2$-bump on $A_k$. 
Moreover, for any
open covering $\{U_i\}$ of $\Omega$ we can in addition insure that each set $B_k$ 
is contained in some $U_i$. 

If $\rho$ does not have any critical points in the set $U_i$ for some $i$
then we can further achieve using Narasimhan's lemma on local convexification 
that for every $k=0,1,\ldots,m-1$ for which $B_k$ lies in $U_i$ it holds that 
$B_k$ is a convex $2$-bump on $A_k$.

Therefore, if $\rho$ does not have any critical points in the set 
$\{z\in \Omega; a\le \rho(z)\le b\}$,
then we can insure that
$B_k$ is a convex $2$-bump on $A_k$ for every  $k=0,1,\ldots,m-1$.
In this case we obtain the map $g$ by using Lemma \ref{K1} $m$ times,
where each time we push the boundary of the disc to the complement of the set $A_j$.

Now assume that $\rho$ has exactly one critical point in the set 
$\{z\in\Omega; a< \rho(z)< b\}$.
Since by assumption this critical point is non-degenerate we can achieve 
that in addition it does not lie on $f(b\disc)$.
Denote this critical point by $q$.
If $q$ is a local minimum then the boundary of the disc already lies above the critical level set; we cannot approach $q$ by the non-critical procedure so we can 
continue in the same way as if there were no critical points.

Therefore from now on we may assume that $q$ is not a local minimum.
We shall push the boundary of the disc to the higher level sets of $\rho$ and we will
keep the boundary away from the critical point. Here we need that the complex 
dimension of the
manifold is at least $3$. In the local coordinates around the critical point
we choose complex $2$ dimensional subspaces $L_1$, $L_2$ and $L_3$ and we move 
the boundary of the disc
in these directions. Away from the critical point we use convex $2$-bumps.
The details are as follows.
We can choose a biholomorphic change of coordinates $\Phi$ from a 
\nbd\ of $q$ to a \nbd\ $\omega$ of $0$ in $\cn$ such that 
$\Phi(q)=0$ and $\rho\circ \Phi^{-1}|\omega$ is strictly \psh\ in the first two coordinates.
Denote by $L$ the complex $2$-dimensional subspace generated by $z_1,z_2$.
Let $\omega_1\ss \omega$ be a \nbd\ of $0$. Then for every small perturbation $L'$ of $L$ and for 
each $z\in \omega_1$ the map $\rho\circ \Phi^{-1}$ is strictly \psh\ on $(z+L')\cap \omega_1$.
Therefore one can choose three complex linear subspaces $L_1$, $L_2$ and $L_3$ in $\cn$ 
such that $L_1\cap L_2\cap L_3=\{0\}$ and for each $j$, $1\le j\le 3$, and for 
$z\in \omega_1$, the map $\rho\circ \Phi^{-1}$ is strictly \psh\ on $(z+L_j)\cap \omega_1$.
There is a $\d>0$ so small that $\{z\in\cn;\dist(z,L_j)<\d,1\le j \le 3\}\ss \omega_1$.
Denote this set by $\omega_0$. By taking smaller $\d$ if necessary we may assume that
$\Phi(f(\bdisc))\cap \omega_0=\emptyset$.

For each $w\in \{z\in \Omega; a\le \rho(z)\le b\}\setminus \Phi^{-1}(\omega_1)$ 
we can choose a coordinate \nbd\ $\Omega_w$, a biholomorphic map 
$\Phi_w\colon \Omega_w\to \Phi_w(\Omega_w)\subset \cn$, where $\Phi_w(\Omega_w)$ is convex 
such that $\rho\circ\Phi_w^{-1}$ is strictly convex in the first two coordinates.
We can further assume
that $\Omega_w\ss \Omega\setminus \bar{\Phi^{-1}(\omega_0)}$.

For each $w\in \Phi^{-1}(\omega_1\setminus \omega_0)$
there is $j$, $1\le j \le 3$, such that $\dist(\Phi(w),L_j)\ge {\d} $. 
By the above  $\rho\circ \Phi^{-1}$ is strictly \psh\ on $(\Phi(w)+L_j)\cap \omega_1$. 
Hence there is a biholomorphic change of coordinates on $\Phi(w)+L_j$ near $\Phi(w)$ such that 
in the new coordinates $\rho\circ \Phi^{-1}$ is strictly convex on $\Phi(w)+L_j$ near $\Phi(w)$.
Since strict convexity is preserved by small perturbations it follows that
there are new coordinates near $\Phi(w)$ in $\cn$ and a \nbd\
$\Omega_w$ of $w$ in $X$ such that $\rho\circ \Phi^{-1}$ is strictly convex in the new coordinates 
on $(\Phi(z)+L_j)\cap \Phi(\Omega_w)$ for each $z\in \Omega_w$. By construction the small tangent discs 
corresponding to $\Phi(z)$
along which we lift in Lemma \ref{K1} lie in $\Phi(z)+L_j$. 
Denote by $\Omega_q$ the set $\Phi^{-1}(\omega_0)$.
Note that $\{\Omega_w\}$ is an open covering of 
$\{z\in \Omega; a\le \rho(z)\le b\}$.
By the above there are a finite number of domains
$\{z\in \Omega;\rho(z)<a\}=A_0\subset A_1\subset\cdots\subset A_m=\{z\in \Omega;\rho(z)<b\}$
such that for each $k$, $0\le k\le m-1$, we have
$A_{k+1}=A_k\cup B_k$, where $B_k$ is a $2$-bump on $A_k$, and 
there is a set $\Omega_w$ such that $B_k\subset \Omega_w$ and if $w\ne q$ 
then $B_k$ is a convex $2$-bump on $A_k$.

We construct the map $g$ inductively.
At each step we construct a continuous map $f_k\colon\cdisc\to X$, \holo\ on $\disc$,
with the following properties 
\begin{enumerate}
\item[(a)] $f_k(\cdisc\setminus r\disc)\cap K=\emptyset$, 
\item[(b)] $f_k(\bdisc)\cap (\bar{ \Omega_q\cup A_k})=\emptyset$, 
\item[(c)] $d(f_k(\z),f(\z))\le\tfrac{k}{2m}\e$ $(|\z|\le r)$,
\item[(d)] $f_k(0)=f(0)$,
\item[(e)] $f_k'(0)=R_kf'(0)$ for some $R_k$, $\root{m}\of{r}<R_k<1$.
\end{enumerate}
Let $f_0=f$ and note that $f_0$ satisfies all the properties. Assume that we have
already constructed the map $f_k$ with the properties (a)-(e) for some
$k$, $0\le k\le m-1$.
If $B_{k}\subset \Omega_q$ then we put $f_{k+1}=f_{k}$. 
In this case the map $f_{k+1}$
obviously satisfies (a), (c), (d) and (e). The property (b) follows from the fact that
$B_{k}\subset \Omega_q$ and that the map $f_k$ satisfies (b).
Otherwise, if $B_{k}\subset \Omega_w$, $w\ne q$, then we use
Lemma \ref{K1} to get the map
$f_{k+1}$. The fact that $f_{k+1}(\bdisc)$ misses $\bar \Omega_q$ follows 
from the properties of the covering;
if $\Omega_w$ misses $\Phi^{-1}(\bar\omega_0)$ and if the perturbation constants
are small enough then obviously $f_{k+1}(\bdisc)$ misses $\bar \Omega_q$.
Otherwise, there is $j$, $1\le j\le 3$, such that 
$\dist(\Phi(w),L_j)\ge {\d} $ and there is a biholomorphic change of coordinates
such that in the new coordinates $\rho\circ \Phi^{-1}$ is strictly convex on 
$\Phi(z)+L_j$ for $z\in\Omega_w$.
Since the boundary of the disc $f_k$ does not intersect
$\bar \Omega_q$, at each point $\z\in b\disc$ such that $f_k(\z)\in\Omega_w$ we have 
$\dist(\Phi(f_k(\z)),L_j)> {\d}$ and
the small tangent disc,
along which we lift the boundary of the disc $f_k$,
lies in $\Phi(f_k(\z))+L_j$.
Therefore, if the perturbation constants
are small enough, it holds that $\dist(\Phi(f_{k+1}(\z)),L_j)> {\d}$ $(\z\in\bdisc)$.
This proves (b).
The properties (a), (c), (d) and (e) are easily satisfied. The construction is finished.
The map $g=f_m$ has all the required properties and the proof is complete.
\end{proof}

\begin{proof}[Proof of Theorem \ref{izrek1}]
By Morse theory (\cite[Observation 4.15]{HL} and \cite[Proposition 0.5]{HL})
we get an exhaustion function $\rho$ of class $\cC^\infty$ without degenerate critical points 
and $M'$ such that the Levi form of $\rho$ has at each point of $\{\rho>M'\}$
at least $2$ positive eigenvalues and such that for $\z\in\bdisc$ it holds
that $\rho(f(\z))>M'$. We may additionally assume that there is only one critical point on
each critical level set. 

Choose an increasing sequence $a_j$ of regular values of $\rho$, converging to $\infty$, and
such that $\rho(f(\z))>a_1$ $(\z\in\bdisc)$ and for each $j\in \N$ there is at most one critical
value on $(a_j,a_{j+1})$. 
Choose a decreasing sequence $\e_j>0$ such that 
\begin{equation}
\text{if }z\in X,\ \rho(z)\le a_j,\ w\in X,\ d(z,w)\le\e_j\text{ then } |\rho(z)-\rho(w)|<1 .
\label{dol}
\end{equation}
Using Lemma \ref{kritdvig} one can construct 
inductively a sequence of continuous maps $f_n\colon\cdisc\to X$, \holo\ on $\disc$,
an increasing sequence $r_n$ of positive numbers converging to $1$,
such that $\sum_{n=1}^\infty (1-r_n)$ converges,
and a sequence $R_n$, $r_n<R_n<1$,
such that for each $n$,
\begin{enumerate}
\item[(a)] $\rho(f_n(\z))>a_n$ $(\z\in\bdisc)$,
\item[(b)] $\rho(f_n(\z))>a_{n-1}-1$ $(r_n\le |\z|\le 1)$,
\item[(c)] $d(f_n(\z),f_{n-1}(\z))<\tfrac{\e_n}{2^n}$ $(|\z|\le r_n)$,
\item[(d)] $f_n(0)=f(0)$,
\item[(e)] $f_n'(0)=\prod_{j=1}^n R_jf'(0)$.
\end{enumerate}
By (c) the sequence $f_n$ converges uniformly on compacta on $\disc$ and 
the limit map $g$ is \holo\ on $\disc$.
Take $n\in\N$ and fix $\zeta$, $r_n\le|\z|\le r_{n+1}$. Then
we get by (c) that 
$$d(f_n(\z),g(\z))\le d(f_n(\z),f_{n+1}(\z))+d(f_{n+1}(\z),f_{n+2}(\z))+\cdots \le $$
$$\le \tfrac{\e_{n+1}}{2^{n+1}}+\tfrac{\e_{n+2}}{2^{n+2}}+\cdots<{\e_n}.$$

Therefore, if $\rho(g(\z))\le a_n$, then this together with (b) and (\ref{dol}) implies that 
$\rho(g(\z))>a_{n-1}-2$.
Since $\lim_{n\to\infty}r_n=1$ and since $\lim_{n\to\infty}a_n=\infty$, it follows that
$g$ is a proper map. 
By (d) we obtain that $g(0)=f(0)$.
Since $\sum_{n=1}^\infty (1-r_n)$ converges, $\sum_{n=1}^\infty (1-R_n)$ converges, and using
\cite[Theorem 15.5]{Rud} we get that the infinite product
$\prod_{j=1}^\infty R_j$ converges to $\l>0$. Therefore $g'(0)=\l f'(0)$.
This completes the proof.
\end{proof}

\textit{Added in the final revision.} 
In the subsequent paper \cite{BDDF} written by Franc Forstneri\v{c} and
the author the conclusion of Theorem \ref{izrek1} is extended 
to complex spaces with singularities.

{\it Acknowledgments.} The author wishes to thank 
F. Forstneri\v{c} for useful
discussions while working on this paper. 
She would also like to thank 
J. Globevnik, M. Slapar, S. Strle and J. Winkelmann for helpful remarks.

This research has been supported in part by 
the Ministry of Education, Science and Sport of Slovenia
through research program Analysis and Geometry, Contract No. P1-0291 and research project
No. J1-6173-0101-04.

\bibliography{biblioAIF1}

\end{document}